\documentclass[letterpaper, 10 pt, conference]{ieeeconf}  

\IEEEoverridecommandlockouts                              

\usepackage{amsmath}
\usepackage{mathtools}
\usepackage{amsfonts}
\usepackage{amssymb}
\usepackage{eurosym}
\usepackage{color}
\usepackage{comment}

\newtheorem{proposition}{Proposition}

\newtheorem{corollary}{Corollary}

\newtheorem{remark}{Remark}

\newcommand{\ba}{\begin{array}}
\newcommand{\ea}{\end{array}}

\newcommand{\be}{\begin{equation}}
\newcommand{\ee}{\end{equation}}

\newcommand{\mc}{\mathcal}

\newcommand{\R}{\mathbb{R}}

\newcommand{\de}{\mathrm{d}}

\newcommand{\tup}[1]{\textup{#1}}

\def\R{\mathbb{R}}
\def\Rp{\mathbb{R}_+}

\usepackage{mathtools}
\usepackage[acronym]{glossaries}
\newacronym{PEV}{PEV}{Plug-in Electric vehicles}
\newacronym{MPC}{MPC}{Model Predictive Control}
\newacronym{SQP}{SQP}{Sequential Quadratic Programming}
\newacronym{MSE}{MSE}{Mean Squared Error}
\newacronym{QP}{QP}{Quadratic Programming}
\newacronym{IBTDM}{IBTDM}{Incentive Based Traffic Demand Management}
\newacronym{FIFO}{FIFO}{first-in-first-out}
\newacronym{TSTT}{TSTT}{total system travel time}

\title{\LARGE \bf
Incentive-Based Electric Vehicle Charging for \\ Managing Bottleneck Congestion
}

\author{Carlo Cenedese$^{1}$, Patrick Stokkink$^{2}$, Nikolas Gerolimins$^{2}$ and John Lygeros$^{1}$ 
\thanks{*This work was supported by the SNSF under NCCR Automation.}
\thanks{
$^{1}$ Carlo Cenedese and John Lygeros are with the Automatic Control Laboratory at ETH Z\:urich, CH-8092 Z\:urich, Switzerland. 
        {\tt\small \{ccenedese, jlygeros\}@ethz.ch}}%
\thanks{$^{2}$ Patrick Stokkink and Nikolas Gerolimins are with Urban Transport Systems Laboratory, School of Architecture, Civil and Environmental Engineering, EPFL, CH-1015 Lausanne, Switzerland 
        {\tt\small \{patrick.stokkink, nikolas.geroliminis\}@epfl.ch}}%
}

\begin{document}

\maketitle
\thispagestyle{empty}
\pagestyle{empty}

\begin{abstract}
We propose an incentive-based traffic demand management policy to alleviate traffic congestion on a road stretch that creates a bottleneck for the commuters. The incentive targets electric vehicles owners by proposing a discount on the energy price they use to charge their vehicles if they are flexible in their departure time. We show that, with a sufficient monetary budget, it is possible to completely eliminate the traffic congestion and we compute the optimal discount. We analyse also the case of limited budget, when the congestion cannot be completely eliminated. We compute analytically the policy minimising the congestion and  estimate the level of inefficiency for different budgets. We corroborate our theoretical findings with numerical simulations that allow us to highlight the power of the proposed method in providing practical advice for the design of policies.
\end{abstract}
\section{Introduction} 
Road traffic is known to be a major source of fuel consumption and a concern for pollution and climate change. The increasing number of road users has lead to an increase in road congestion over the years, and as a consequence the concentrations of several pollutants are expected to nearly double during the rush hour in modern cities \cite{zhang2011vehicle}. The fast share  of \glspl{PEV} is also negatively affected by  congestion, in fact they experience a higher battery consumption~\cite{bingham2012impact}. 

Highly congested periods are usually associated with commuting times during the morning and evening. 
The bottleneck model has been developed to capture situations in which  commuters face one main source of road congestion during their travel. Arguably, the most relevant application can be found in the modelling of a highway that connects origins to destinations \cite{vickrey1969congestion}. The structural model explicitly incorporates physical aspects of the road congestion as well as behavioural decisions of drivers. Individual commuters select their departure time as a best response to their local travelling cost that depends on the  departure time chosen by the others. Due to this selfish inclination of the drivers, the road capacity can be  saturated leading to congestion \cite{arnott1990economics, arnott1993structural,lindsey:2004:existence_uniqueness_of_bottleneck_model_equilibrium, small2015bottleneck}. We refer to~\cite{li2020fifty} and references therein for a recent overview on applications and findings regarding the bottleneck model.

In the literature, various solutions have been proposed to reduce bottleneck congestion. On the one hand, there are policies that enforce a price on the use of the bottleneck segment, e.g.,  tolling \cite{van2011congestion, vandenberg:2011:winning_or_losing_from_bottleneck, lindsey2012step}. Such policies are known as \textit{hard policies} as the users are forced to comply.   On the other hand, there are \textit{soft policies}, also known as \gls{IBTDM},  where commuters are rewarded if they enter the bottleneck outside the peak of congestion, \cite{sun2020managing}. Whereas policies such as tolling have been extensively studied, positive incentives have received considerably less attention. Among others, the question of how incentives can be efficiently implemented in real-life situations remains open. 

\smallskip
In this work, we propose an  \gls{IBTDM} policy where incentives are offered to \gls{PEV} owners if they are flexible in their departure time and their associated cost of early or late arrival at the destination. The reward is provided via a discount on the electricity price. 
We show that such a policy is effective in reducing the \gls{TSTT} and derive the total monetary budget necessary to implement it. We also analytically compute the level of inefficiency that arises due to a diminishing return of the monetary budget, connected to the commuters' perception of the value of the incentive.   Moreover, we study the case in which only a  limited budget is allocated to implement the \gls{IBTDM},  in which case the congestion cannot be completely eliminated. Also in this case, we are able to derive the optimal resource allocation and the associated policy.

Notice that by offering incentives to \gls{PEV} owners, our policy is targeted on a specific type of road users. This approach is common in the literature, for example in carpooling. In the presence of High Occupancy Vehicle (HOV) lanes, toll-differentiation can be applied to reduce congestion \cite{yang1999carpooling}.
\smallskip

The remainder of the paper is structured as follows. In Section~\ref{sec:prob_form}, we introduce briefly the classical bottleneck model and we extend it to comprise the charging of \glspl{PEV}. In Section~\ref{sec:optimal_policy}, we derive the optimal discount for the energy price in the case of unlimited and limited budget. Moreover, we derive the monetary budget and the associated inefficiency. Section~\ref{sec:num_sim} provides numerical simulations that corroborate our theoretical findings. We conclude the paper by discussing possible extensions of the proposed methodology.

\section{Problem formulation}
\label{sec:prob_form}
\subsection{Classical bottleneck model}
We briefly review the classical bottleneck model that is the foundation of our analysis. The discussion follows \cite{arnott:1987:schedule_delay_departure}.  The model comprises $N\in\mathcal{N}$ commuters (or agents) who travel from their origin, e.g., home, to their destination, e.g., work. During their trip, they pass through a single bottleneck, that is assumed to be the only potential source of congestion they may encounter. The capacity of the road at the bottleneck is assumed to be constant and equal to $s>0$ vehicles per time instant.
If the departure rate $r(t): \R_+\rightarrow\R_+$ at which the vehicles enter the bottleneck is greater than $s$,  for some time instant $t$, then a queue $Q(t) : \R_+\rightarrow \Rp$ is created. The queue dynamics read    	 \smallskip
\begin{equation}\label{eq:queue_integral}
Q(t) = \int_{\hat t}^t r(\tau)\, \de\tau - s(t-\hat t)\:,
\end{equation}
\smallskip
where $\hat t <t$ is the last moment at which there was no congestion, i.e., $Q(\hat t)=0$.  The vehicles leave the queue according to a \gls{FIFO} principle. For each agent $i\in\{1,\dots, N\} = \mathcal{N}$ arriving at the bottleneck 
 at time $t$, the complete travel time experienced is
\begin{equation}
\label{eq:travel_time}
T_i(t) = T^{\tup f}_i+T^{\tup v}(t),
\end{equation}
where $T^{\tup f}_i\geq0$ denotes the fixed time it takes agent $i$ to commute in the absence of traffic congestion, while $T^{\tup v}(t)\geq 0$  is the additional time spent due to the traffic congestion; clearly $T^{\tup v}(t)> 0$ if and only if $Q(t)>0$.  As every agent can have a different origin and/or destination, each agent can have a different commuting path and therefore  $T^{\tup f}_i$ can vary for different $i$. On the contrary, the value of $T^{\tup v}(t)$ depends only on the time at which the vehicles enter the bottleneck. 
The time $T_i^{\tup f}$ does not play a role in the agents' decision-making process, so without loss of generality we assume that $T^{\tup f}_i=0$ for all $i\in \mc N$, as done in \cite{arnott:1987:schedule_delay_departure}.  
The waiting time at the bottleneck is 
\begin{equation}\label{eq:travel_time_w_Q}
T^{\tup v}(t) = \dfrac{Q(t)}{s} .
\end{equation}
Therefore, if we denote by $t^*$ the common desired arrival time of the commuter, as in~\cite{arnott:1987:schedule_delay_departure}, the time \mbox{$0\leq t'\leq t^*\leq N/s$} at which it has to enter the bottleneck (or equivalently departure since $T^{\tup f}_i=0$) to reach its destination at $t^*$ is 
\begin{equation}
t' + T^{\tup v}(t') = t^*.
\end{equation}
 
The agents are assumed to be perfectly rational, that is, they select the $t$ that minimises their own ``cost''. The cost of an agent is associated to the discomfort experienced by choosing a particular departure time.
If a commuter leaves at $t<t'$ or $t>t'$, then it will arrive early or late, respectively. The cost per unit of time of arriving early is $\beta>0$, that of being late is $\gamma>0$, while time spent in congestion is penalized by $\alpha > 0$ \cite{vickrey1969congestion,arnott:1987:schedule_delay_departure}. 
Therefore, the cost that each commuter aims to minimise is 
\smallskip
 \begin{equation}\label{eq:cost_i}
 \begin{split}
 \hat C(t) =& \alpha T^{\tup v}(t) +   \max\left(\beta d(t,t^*),-\gamma  d(t,t^*) \right),\\
d(t,t^*) \coloneqq & t^*-t-T^{\tup  v}(t).
 \end{split}\end{equation}
As shown in \cite{arnott:1987:schedule_delay_departure,vandenberg:2011:winning_or_losing_from_bottleneck}, one should select $\beta< \alpha < \gamma$ to avoid unnatural user behaviour.
Note that we consider homogeneous agents; we comment on the case where commuters are heterogeneous in Section~\ref{sec:conclusion}.


\subsection{Extended bottleneck model with charging of \glspl{PEV}}
We extend the classical bottleneck model by considering the presence of \glspl{PEV} that receive an incentive if they do not travel during congested periods. The incentive is a discount $p(t):\R_+\rightarrow \R_+$ over the price of the electricity purchased by the owner. Here, $p(t)$ depends on the time at which the \gls{PEV} enters the bottleneck, rather than the time at which it starts charging. We assume that each \gls{PEV} must  charge for $\bar\delta>0$  time instants before entering the bottleneck. They can choose to charge at home, at a fixed energy price $\bar p$, or, at one of the charging stations that take part in the policy at a variable energy price $\bar p -p(t) $. The charging station can be located anywhere along the path that the commuter takes to reach the bottleneck. This addition endows the \glspl{PEV} commuting with an additional degree of freedom. In fact, they  now have two decision variables: 
\begin{enumerate}
    \item the time at which they enter the bottleneck $t$,
    \item the time spent at the charging station  \mbox{$\delta\in[0,\bar \delta]$.}
\end{enumerate} 
As a consequence of the above definitions, the charging station is entered at time $t - \delta$ and the remaining electricity is assumed to be received at home, see Remark~\ref{rem:chargign_home} for an alternative formulation. 
 The price discount $p(t)$ can be implemented through a booking system, where a charging spot for the time period $[t - \delta,t]$ is reserved in advance. 

\begin{figure}[t]
    \centering
    \includegraphics[width = 0.5\textwidth]{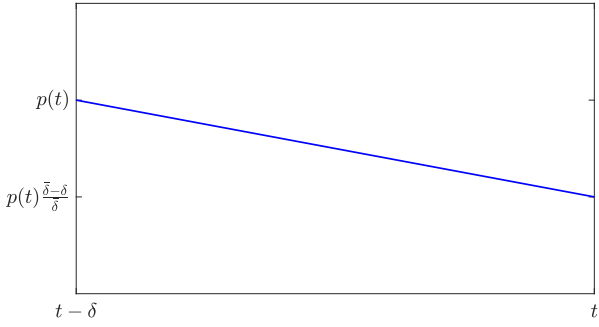}
    \caption{Perceived incentive on the interval $[t - \delta, t$], where charging starts at $t - \delta$ and ends at $t$.}
    \label{fig:perceivedincentive}
\end{figure}

Even though the discount is constant, we assume that commuters perceive a diminishing benefit the longer they charge. This is due to the increasing battery level that leads to a lower  necessity to charge. Furthermore, it also models the well-known psychological phenomenon  of the {\it range anxiety} that affects \gls{PEV}'s owners, whereby they tend to experience a growing discomfort the more the battery of their vehicle discharges, see \cite{thomas:2011:range_anxiety,rauh:2015:range_anxiety_2}. Consequently, the longer the commuters  charge at the charging station, the less appealing the incentive becomes. Since the battery dynamics is often modelled as an integrator with respect to the purchased energy, we assume that this reduction is  linear over $[t - \delta,t]$. If $t$ is the time at which an agent leaves the charging station and enters the bottleneck and, $\delta$ the charging time, then the perceived incentive at time $\tau \in [t - \delta, t]$ is defined as 
\begin{equation}
\label{eq:p_hat}
   \hat{p}(\tau, t) \coloneqq - \frac{p(t)}{\bar{\delta}}\tau + \frac{p(t)}{\bar{\delta}}(t - \delta + \bar{\delta}).
\end{equation}
As shown in Figure~\ref{fig:perceivedincentive}, the perceived discount at the end of the charging performed at the charging station is $p(t)(\bar \delta-\delta)/\bar \delta<p(t)$.
The time $t$  depends both on the potential incentive, as well as, on the delay, earliness and lateness imposed by the bottleneck. Notice that the linearity of the diminishing benefit plays also a technical role, since it allows us to derive closed form solutions for most of the optimal quantities computed in the remainder.

Next, the cost in~\eqref{eq:cost_i} is modified to take into consideration these new features of the model, leading to 
\begin{equation}\label{eq:cost_t_delta}
C(t, \delta) = \hat C(t) + C_{\tup{ch}}(t, \delta),
\end{equation}
where
\begin{equation}
\label{eq:cost_charge}
   C_{\tup{ch}}(t, \delta)  \coloneqq \alpha \delta + (\bar \delta-\delta)\bar p + \int_{t - \delta}^t (\bar p -\hat{p}(\tau,t)) d\tau,  
\end{equation}
 which comprises the inconvenience of waiting at the charging station for $\delta$ time units, added to the perceived incentive $\hat p$. 

The optimal charging time $\delta^*(t)$  depends on the bottleneck through the departure time $t$. We can therefore determine the optimal charging time $\delta^*(t)$ that minimises the charging cost for the agents departing at time $t$. To this end, we substitute~\eqref{eq:p_hat} and \eqref{eq:cost_charge} in~\eqref{eq:cost_t_delta}, and then optimise $C$ with respect to its second argument $\delta$ over the interval $[0,\bar\delta]$. The resulting optimal charging time reads
\begin{equation}
\label{eq:delta_star}
    \delta^*(t) = \max\left\{\left(1 - \frac{\alpha}{p(t)}\right) \bar{\delta},0\right\}. 
\end{equation}
It can be easily verified that $\delta^* < \bar{\delta}$ if $ \frac{\alpha}{p(t)} > 0 $. In case $p(t) < \alpha$, $\delta^* = 0$ as the value of time spent waiting is higher than the received discount.
Substituting in~\eqref{eq:cost_charge}, we recast the cost  as a function of $p(t)$ and $t$, as 
\begin{equation}
\label{eq:cost_charge_end}
   C_{\tup{ch}}(t) = \begin{cases}
   -\frac{\bar \delta}{2p(t)} (p(t)-\alpha)^2 + \bar \delta\bar p ,& \:\text{if} \:p(t)\geq \alpha\\
   \bar \delta\bar p ,& \:\text{otherwise},
   \end{cases}
\end{equation}
where we have dropped the second argument since it is  always $\delta=\delta^*$.
\smallskip 
\begin{remark}[Constant perceived incentive]
\label{rem:no_dim_ret}
If the provided incentive $p(t)$ does not suffer from a diminishing return, as in \eqref{eq:p_hat}, then the cost of charging simplifies to
\begin{equation}
\label{eq:cost_charge_tmp}
   \hat C_{\tup{ch}}(t, \delta)  \coloneqq (\alpha-p(t)) \delta + \bar \delta\bar p   ,  
\end{equation}
that is linear in $\delta$. In this case, solving for $\delta^*$ is trivial: $\delta^*=0$ if $p(t)<\alpha$, and $\delta^*=\bar \delta$ whenever $p(t)>\alpha$. Notice that the analysis carried out in the remainder of the paper can be easily modified and simplified to address this degenerate case. \hfill\QEDopen
\smallskip
\end{remark}
\begin{remark}[Equivalent problem setup]
\label{rem:chargign_home}
   In \eqref{eq:cost_charge_end} the constant component $\bar \delta \bar p$  does not play a role in the choice of the  optimal pricing strategy, since it is experienced by the \glspl{PEV} independently from whether or not they stop at the charging station. It is interesting to notice that the problem in which the users can charge for at most $\bar \delta $ time intervals, but do not charge at home,  shares the same optimal pricing with the setup described above, both in the case of limited and unlimited budget. In this alternative problem formulation, the cost of charging reads as
   $$\tilde C_{\tup{ch}}(t, \delta)  \coloneqq \alpha \delta + \int_{t - \delta}^t (\bar p -\hat{p}(\tau,t)) d\tau\,.$$
   This interpretation suites better the case in which the policy maker wants to provide the incentive independently from where the commuter charges the vehicle, e.g., at the charging station or at home. On the contrary, the proposed implementation considers the charging stations as the main actor in providing the incentive.  
\hfill\QEDopen
\end{remark}

\section{Optimal Pricing Policy}
\label{sec:optimal_policy}
\subsection{Unlimited Budget}
To compare pricing policies, we introduce the concepts of \gls{TSTT} and budget. The \gls{TSTT} reflects the time that commuters cumulatively spend in congestion and is defined as $\Gamma \coloneqq \int_0^{N/s} T^{\tup v}(t)\, dt$. The budget is the total amount of money earmarked by the policy maker to provide the discount $p(t)$ on the energy price and is computed as 
\begin{equation}
\label{eq:M_euro}
M_\$ = \int_0^{\frac{N}{s}} r(\tau) \delta(\tau)p(\tau) d\tau\:.
\end{equation} 
The diminishing effect of the incentive throughout the charging period $\delta(t)$ creates some inefficiency in the allocation of $M_\$$. Loosely speaking, the effect that a given incentive has in affecting the commuters' behavior diminishes the longer they stop at the charging station. We denote by  $M_{\tup {per}}$ the part of $M_\$$ that is perceived by the users and define this \textit{perceived incentive} as 
\begin{equation}
\label{eq:M_p}
M_{\tup {per}} = -\int_0^{\frac{N}{s}} r(\tau) (C_{\tup{ch}}(\tau)-\bar\delta\bar p) d\tau\:,
\end{equation}
where, following \eqref{eq:cost_charge_end}, $C_{\tup{ch}}(t)-\bar\delta\bar p$ is the incentive perceived by a \gls{PEV} entering the bottleneck at $t$ and stopping at the charging station for $\delta^*(t)$ instants. Notice that this definition intrinsically assumes  that the commuters adopt the policy $\delta^*$ while that in \eqref{eq:M_euro} does not.
By substituting \eqref{eq:cost_charge_end}~and~\eqref{eq:M_euro} in~\eqref{eq:M_p}, we can derive the following  relation between the two budgets 
\begin{equation}
\label{eq:M_per}
        M_{\tup {per}}  =M_\$  -\int_0^{\frac{N}{s}}r(\tau)\left[ \frac{p(\tau)\delta^*(\tau)^2}{2\bar \delta} +\alpha\delta^*(\tau)  \right]\, d\tau.
\end{equation}
The relation can be refined further by substituting   \eqref{eq:delta_star} to obtain
\begin{equation}
\label{eq:Delta_M1}
\begin{split}
       \Delta_M&\coloneqq  M_\$ - M_{\tup {per}}\\ 
        &= \int_0^{\frac{N}{s}} r(\tau) \max\left\{\frac{p(\tau)^2-\alpha^2}{2p(\tau)},0\right\}\, d\tau \geq 0\,,
        \end{split}
\end{equation}
where the inequality holds since  $r(t)\geq 0$ for all $t$, and the maximum is due to \eqref{eq:delta_star} that implies that for $p(t)<\alpha$ no agent perceives the incentive, and thus $M_\$=M_{\tup{per}}=0$. The above quantity describes the level of inefficiency in how $M_\$$ influences the users' behaviour. We call this the \textit{inefficiency gap} and denote it by $\Delta_M\geq0$. Notice that $\Delta_M=0$ only in degenerate cases, namely if no commuters travel ($r(t)=0$) or no one charges at the charging station ($\delta^*=0$).

A policy $p(t)$ is defined to be optimal if it minimises $\Gamma$ while using the smallest budget $M_\$$. 
The following proposition provides the analytical formulation of the optimal policy, depicted in Figure~\ref{fig:price}, and of the associated budget.
\smallskip
\begin{proposition}[Optimal discount with unlimited budget]
\label{prop:opt_policy}
Consider the extended bottleneck model  introduced above, then the optimal energy price discount is 
\begin{align}
\label{eq:p_1}
    p^*(t) = \alpha - \frac{g}{\bar{\delta}}  + \sqrt{\left(\frac{g}{\bar{\delta}}-\alpha\right)^2 -  \alpha ^2} \\\nonumber 
     \text{where } g\coloneqq \begin{cases}
     \beta (t-t^*),\, \text{if } t<t^*\\
     \gamma ( t^*-t),\, \text{otherwise}
     \end{cases}\,,
\end{align} 
the \gls{TSTT} is $\Gamma=0$, and the required budget 
\begin{equation}\label{eq:M_eu_star}
    M_{\$}^* =  \frac{\beta\gamma N^2}{2s(\beta+\gamma)} + s \int_0^{\frac{N}{s}} \frac{p^*(\tau)^2-\alpha^2}{2p^*(\tau)} \,d\tau.
\end{equation} 
\end{proposition}
\smallskip
\begin{proof}
See Appendix \ref{app:appendix_proof_1}.
\end{proof}
\smallskip
As can be seen in Figure~\ref{fig:price} and in accordance with \cite{sun2020managing}, the maximum  incentive is received the further $t$ is from $t^*$, while it reaches its minimum at $t^*$, namely $p^*(t^*)=\alpha$ that implies $\delta^*(t^*)=0$ from \eqref{eq:delta_star}. Moreover, the fact that at the equilibrium the queuing delay is eliminated by $p^*(t)$ implies that $r(t)=s$ for all $t\in[0,N/s]$.
By applying the optimal policy $p^*$, we can specify the formulation of the efficiency gap in \eqref{eq:Delta_M1} as
By substituting $p^*$ into \eqref{eq:Delta_M1}, we obtain the efficiency gap at optimality
\smallskip
\begin{align}
\label{eq:Delta_M_star}
    \Delta^*_M &=  s \int_0^{\frac{N}{s}} \frac{p^*(\tau)^2-\alpha^2}{2p^*(\tau)}\,d\tau .
\end{align}
The value above together with $M_\$^*$ provide an important estimation of the resources necessary to eliminate  congestion. In particular, they allow the policy maker to compute the efficiency of the policy, and to estimate the  monetary burden to eliminate the  congestion. It is important to stress that the integrals in \eqref{eq:M_eu_star} and \eqref{eq:Delta_M_star} admit a closed form solution obtained via calculations  similar to those in Appendix~\ref{app:deriv_f}.

\subsection{Limited Budget}
In most real-world scenarios, the monetary budget that is allocated for implementing the proposed incentive will be limited by external factors. If the available budget is lower than $M^*_\$$, then  it will not be possible to completely eliminate congestion. Here we derive a pricing policy that minimizes the \gls{TSTT} subject to such budget constraints.

The following proposition provides the optimal incentive $p(t)$, and the associated departure rate $r(t)$ attained at the equilibrium when a limited budget $M_\$\leq M^*_\$$ is provided. We use $M_\$ = f(M_{\tup{per}}) $ where the definition and the derivation of $f$ is in Appendix~\ref{app:deriv_f}. Even though such relation cannot be directly inverted, one can easily compute the value of $M_{\tup{per}}$ associated to a certain $M_\$$ numerically. Finally,  we use here $t''$, instead of $t'$ (recall that $t'$ is the departure time to arrive exactly at $t^*$, if no \gls{IBTDM} is implemented), to highlight that this is the same quantity in a different setting.   
\smallskip

\begin{proposition}[Optimal discount with limited budget]\label{prop:opt_disc_lb}
Given a monetary budget $M_\$$, the optimal incentive for $t\in[0,N/s]$ is 
\begin{align}
    \label{eq:p_limited_budget}
    & p(t) = \begin{cases} 0 &  t\in(t^\ell,t^r)\\
         \alpha -\frac{g}{\bar \delta}+\sqrt{\left(\frac{g}{\bar \delta}-\alpha\right)^2-\alpha^2}& \text{otherwise }
    \end{cases}\,,\\\nonumber
&  
\text{with } g\coloneqq \begin{cases}
     \beta (t-t^\ell),\, \text{if } t\leq t^\ell\\
     \gamma ( t^r-t),\, \text{if } t\geq t^r
     \end{cases}
\end{align}
and
\begin{equation}\label{eq:t_l_and_t_r}
\begin{split}
t^\ell &= \sqrt{\frac{2\gamma f^{-1}(M_{\$})}{s\beta(\beta+\gamma)}} \\
t^r &= \frac{N}{s}-\sqrt{\frac{2\beta f^{-1}(M_{\$})}{s\gamma(\beta+\gamma)}} \,.
\end{split}
\end{equation}
The associated departure rate is
\begin{equation}
\label{eq:departure_time_lb}
    r(t)=\begin{cases}
    s\,, &  t\in[0,t^\ell]\cup [t^r,N/s]\\
    \alpha s/(\alpha-\beta)\,, &  t\in(t^\ell,t'']\\
    \alpha s/(\alpha+\gamma)\,, &  t\in(t'',t^r)
    \end{cases}
\end{equation}
with $t''=t'+\frac{1}{\alpha}\sqrt{\frac{2\beta\gamma f^{-1}(M_{\$ })}{s(\beta+\gamma)}}$.
\hfill
\end{proposition}
\smallskip
\begin{proof}
See Appendix \ref{subsec:appendix_proof_2}.
\end{proof}
As expected, the incentive is positive only outside the central part of the rush hour. Consequently those periods are the only ones during which there is no congestion and the departure rate coincides with the road capacity, see \eqref{eq:departure_time_lb}. The discontinuity is an effect of the optimal charging time in Equation \eqref{eq:delta_star}, that suggests that any discount lower than $\alpha$ will not convince the commuters to charge. The queue starts to form for $t>t^\ell$, it peaks at $t=t''$ and it  disappears after $t^r$. The value of $M_\$$ grows with $M_{\tup{per}}$, thus the region in which no incentive is provided shrinks the higher the budget. If  $f^{-1}(M_{\$ })=M_{\tup{per}}^*$, then we retrieve the results attained in the previous section, i.e., $t^\ell=t^r=t^*$.
The result above allows us to derive the \gls{TSTT} in the case of limited budget.

\smallskip
\begin{corollary}[\gls{TSTT} with limited budget]
Given a budget $M_\$<M_\$^*$, if the incentive policy $p(t)$ is as in \eqref{eq:p_limited_budget}, then the \gls{TSTT} becomes  
\begin{equation}
    \label{eq:Gamma_lb}
    \Gamma=\frac{1}{\alpha} f^{-1}(M_{\$})-\theta \sqrt{f^{-1}(M_{\$})}+\nu>0\,,
\end{equation}
where \mbox{$\theta\coloneqq (N/\alpha)\sqrt{2\beta\gamma/[s(\beta+\gamma)]}$} and \mbox{$\nu\coloneqq \beta\gamma N^2/(2s\alpha(\beta+\gamma))$}.
\hfill
\end{corollary}
\smallskip

\begin{proof}
From the definition of $\Gamma$, we  notice that that
 \begin{equation}
     \label{eq:gamma_lim}
\Gamma= \frac{s}{2}(t^r-t^\ell)(t^*-t'')\,.
 \end{equation}
So, the proof is completed by substituting \eqref{eq:t_l_and_t_r} into  \eqref{eq:gamma_lim}.
\end{proof}\smallskip
Following \cite[Sec.~3.2]{sun2020managing}, and  taking the first and second derivative of $\Gamma$ with respect to $M_{\tup{per}}$, one can show that the value of $\Gamma$ decreases as $M_{\tup{per}}$ increases, but at a decreasing rate. This implies that smaller values of  $M_{\tup{per}}$, and in turn also of $M_\$$, create greater gains in welfare per budget unit.


\section{Numerical Results}
\label{sec:num_sim}
To quantify the effect of implementing the IBTDM policy for PEV owners, we consider a numerical example. We use the same parameter values as \cite{arnott1990economics} and \cite{lindsey2012step}, that were obtained based on the analysis of \cite{small1982scheduling}. The unit cost parameters are $\alpha = 6.4[\$/\text{h}]$, $\beta = 3.9[\$/\text{h}]$ and $\gamma = 15.21[\$/\text{h}]$, the number of commuters $N = 9000$, and the bottleneck capacity $s = 60 [\text{veh.}/\text{min}]$; thus the period considered lasts $T=150$ min. We assume that the \glspl{PEV} must charge for $\bar \delta=20$ min.\\
\begin{figure}[t]
    \centering
    \includegraphics[width = \columnwidth]{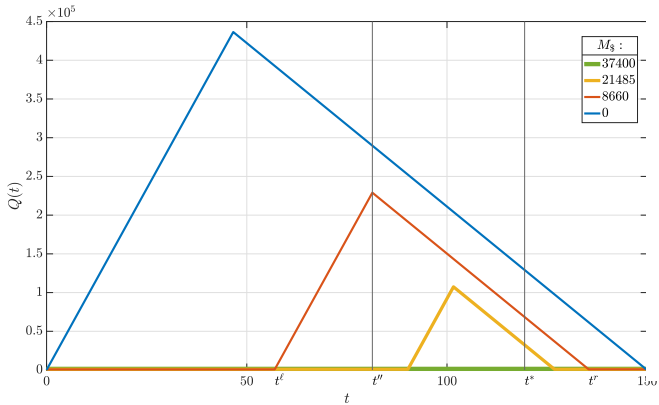}
   \caption{The evolution of the queue $Q(t)$ at the bottleneck for growing values of the budget $M_\$$, where $M_\$^*=\$37400$. We indicate the starting \mbox{$t^\ell=57 \,\text{min}$} , ending $t^r=135.2\,\text{min}$ and peak $t''=81.3\,\text{min}$ instants  of the queue for \mbox{ $M_\$=\$ 8660$}.}
    \label{fig:queue}
\end{figure}\begin{figure}[t]
    \centering
    \includegraphics[width = \columnwidth]{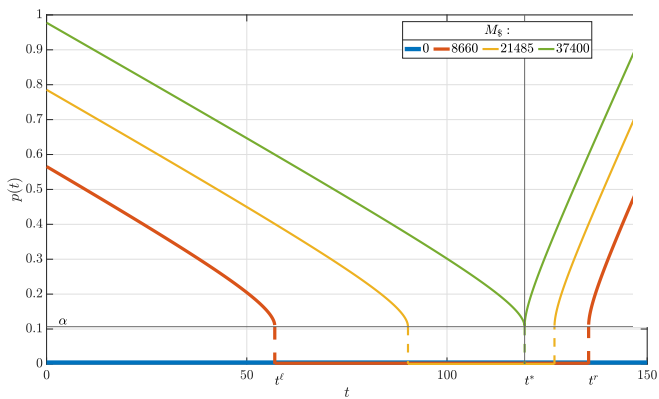}
    \caption{The optimal  discount $p(t)$ for the electricity price for growing values of the budget $M_\$$, where $M_\$^*=\$37400$. For $M_\$=\$ 8660$,  we again indicate the starting $t^\ell=57 \,\text{min}$ and ending $t^r=135.2\,\text{min}$ instants during which $p(t)=0$.}
    \label{fig:price}
\end{figure}
\indent In Figure~\ref{fig:queue}, the evolution of the queue over $[0,T]$ is  depicted for different budgets. When there is no incentive ($M_\$ =0$) the queue peaks at $Q(t)=4.3\cdot 10^{5}$. As shown in \cite{arnott1993structural}, the maximum is reached at $t=t'$. As anticipated, the implementation of the \gls{IBTDM} reduces the congested period $[t^\ell,t^r]$. For increasing values of  $M_\$$, the time $t''$ at which $Q(t)$ peaks move closer to $t^*$, as foreseen in Proposition~\ref{prop:opt_disc_lb}. Moreover, the magnitude of the peak decreases, until congestion disappears for $M_\$^*=\$37400$. Note that even if the available budget is $M_\$ = \$8660$ (less than $25\%$ of $M_\$^*$), the congested period is reduced from $150$ min to $t^r-t^\ell=78$ min (around $48\%$). Furthermore, the maximum length of the queue almost halves to $2.28\cdot 10^{5}$ in this scenario. 

In Figure~\ref{fig:price}, the optimal discount $p(t)$ from Propositions~\ref{prop:opt_policy}~and~\ref{prop:opt_disc_lb} is plotted for the same values of $M_\$$ used for Figure~\ref{fig:queue}. The incentive $p(t)$ is zero during $[t^\ell,t^r]$ and positive (indeed at least $\alpha$) otherwise. In the case of sufficient budget, the incentive becomes  $p(t^*)=\alpha$ and thus $\delta^*(t^*)=0$, see \eqref{eq:delta_star}.\\    
\begin{figure}[t]
    \centering
    \includegraphics[width = \columnwidth]{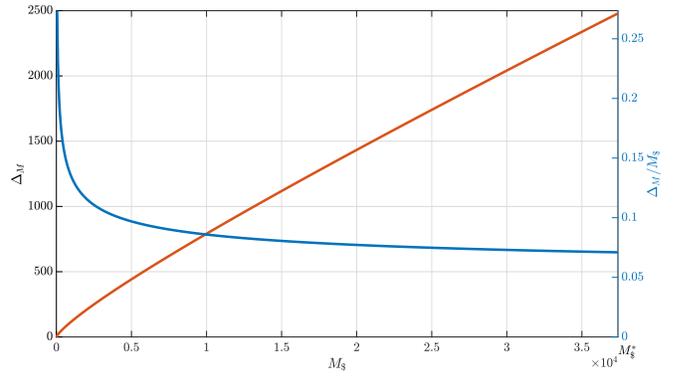}
    \caption{The variation of the inefficiency gap $\Delta_M$ (in blue) and $\Delta_M/M_\$$ (in red) with respect to $M_\$$. }
    \label{fig:delta_M}
\end{figure}\begin{figure}[t]
    \centering
    \includegraphics[width = \columnwidth]{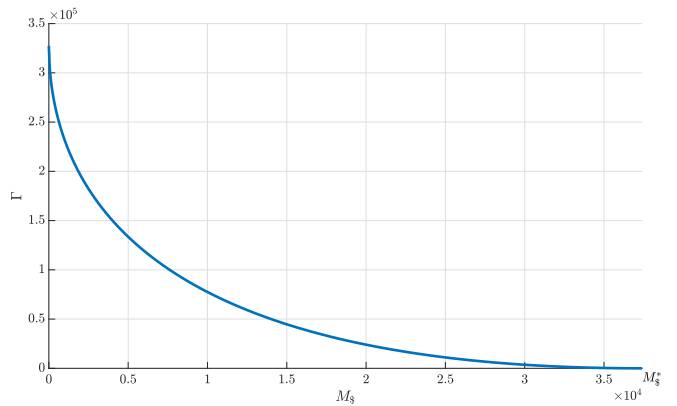}
    \caption{The \gls{TSTT} evolution with respect to different values of monetary budget $M_\$$.}
    \label{fig:gamma}
\end{figure}
\indent Insights on how to efficiently select the budget $M_\$$  can be inferred from Figure~\ref{fig:delta_M}, that shows the relation between the budget and the (relative) inefficiency gap. Even though in absolute terms the inefficiency increases as the budget increases, it decreases as a percentage of the budget. Moreover, the reduction in percentage inefficiency is steeper for smaller values of $M_\$$ while the curve flattens the closer $M_\$$ is to $M_\$^*$. As a consequence, a value of $M_\$$ that is too small generates high inefficiency and therefore should be avoided. 

Figure~\ref{fig:gamma} shows that the \gls{TSTT} monotonically decreases with respect to $M_{\$}$. For $M_\$ = \$8660$, we achieve a remarkable reduction of roughly $70\%$ of its value, see Figure~\ref{fig:gamma}. This should not be surprising since $\Gamma$ is related to the area under $Q(t)$ which decreases rapidly with $M_\$$, see Figure \ref{fig:queue}. It is worth noting that, by considering  the percentage of $\Delta_M$ in Figure~\ref{fig:delta_M}, we can notice a knee around $M_\$=\$ 5000$ after which the reduction in inefficiency percentage is not as fast. Higher social welfare and lower percentage inefficiency are attained for higher budgets, therefore the policy maker should avoid the use of very low budget since they might lead to unsatisfactory results.  The actual value of the best budget is not necessarily obvious a priory, but it depends on the particular design objective and the relative importance of congestion reduction and money spent.

\section{Conclusions}
\label{sec:conclusion}
The use of a dynamically discounted energy price as an incentive for commuters that own a \gls{PEV} is an effective policy to decrease bottleneck congestion during rush hours. In the case of sufficient monetary budget, a policy maker can completely eliminate congestion during the rush hour. Otherwise, the optimal policy, that is provided only during the tails of the original congested period, minimises the duration and the intensity of the traffic jam. The efficiency of the policy can be analytically quantified. It depends on the diminishing return the commuters perceive when they increase the period during which they receive the incentive.

By building on the results on \cite{sun2020managing}, one can extend the findings presented here to the case of heterogeneous groups of commuters, i.e., with different values of time and costs for arriving late/early. Another variation of the proposed approach is considering a limited penetration rate of the proposed policy. Intuitively, there is a strong bond between this case and the limited incentive case. In fact, in the latter some users do not perceive any discount, and thus they can be interpreted as non-compliant. A detailed analysis of such a topic is left for future works. Finally, it is of high interest  to combine the proposed approach with other traffic control mechanisms. For example, the use of tolling can make the proposed policy self-sustainable. Alternatively, it is compelling to study the case in which, on top of the monetary incentive, the vehicles have also the possibility to travel via priority lanes. The arising question in this framework would be the optimal splitting of the road capacity given the incentive already provided.

\section{Appendix}

\subsection{Proof of Proposition~\ref{prop:opt_policy}}\label{app:appendix_proof_1}
From \eqref{eq:cost_t_delta} and \eqref{eq:delta_star}, the total travel cost for each agent is $C(t)=\hat C(t)+C_{\tup{ch}}(t)$, where we obtained  $C_{\tup{ch}}(t)$ as $C_{\tup{ch}}(t,\delta^*)$. Notice that $C_{\tup{ch}}(t,\delta^*)-\bar\delta\bar p$ represents the totality of the incentives and costs due to the policy put in place, since $\bar \delta \bar p$ is experienced by the drivers independently form whether or not they take advantage of the discount or not.  
The proof is performed in two steps, we first derive a form of $C_{\tup{ch}}(t)$ that leads to $\Gamma=0$ and uses the least perceived budget, and then we derive the associated $p^*(t)$.

We consider a linear form
\begin{align}
\label{eq:incentive_p1}
-\psi(t)\coloneqq C_{\tup{ch}}(t)-\bar \delta \bar p = 
    \begin{cases}
    \beta t - g_1 &t \in [0, t'] \\
    -\gamma t - g_2 &t \in (t', N/s].
    \end{cases}
\end{align}
The continuity and differentiability of $C(t)$ with respect to $t$ in the intervals $[0,t')$ and $(t',N/s]$ is guaranteed by \cite[App.~B-C]{sun2020managing}.
Then, by the stationary condition of the optimum, i.e.,$ \frac{dC(t)}{dt}=0 $ 
and substituting \eqref{eq:queue_integral} and \eqref{eq:travel_time_w_Q}, we can show that $r(t)=s$ for all $t\in[0,N/s]$.  It follows that $\Gamma=0$ at the optimum.

For the proposed policy to be an incentive (and not a cost), it must hold that $\psi(t) \geq 0$ for all $t$. Since $\Gamma=0$, the minimum is reached at $t=t^*$. From \eqref{eq:M_per}, it is follows that the perceived incentive used for the implementation, i.e.,  $M_{\tup{per}}^*$, is minimum if $\min_t\psi(t)=0$, see \cite[Prop.~1]{sun2020managing}. Consequently, we obtain that $g_1=\beta t^*$ and $g_2 = -\gamma t^*$. Thus this incentive is socially optimal and uses the least perceived budget.

Next, we retrieve the optimal price discount $p^*(t)$ by substituting \eqref{eq:cost_charge_end} into \eqref{eq:incentive_p1}. For $t\in[0,t')$, we obtain the equation
\begin{equation}\label{eq:equate_c_charge_and_opt}
   \max\left\{ \frac{\bar \delta}{2p(t)} (p(t)-\alpha)^2,0\right\} = \beta (t -  t^*).
\end{equation}
The price discount should always be positive $p(t)\geq 0$ and greater than the value of time, i.e., $p(t)>\alpha$. By solving the relation above for $p(t)$, we obtain the optimal discount policy for $t\in[0,t^*)$ associated to an incentive in the form of $\psi(t)$, so

\begin{align}
\label{eq:p_11}
    p^*(t) = \alpha - \frac{1}{\bar{\delta}} \beta (t -  t^*) + \sqrt{b^2 -  \alpha ^2}, 
\end{align} 
where $b\coloneqq -\alpha +  \frac{1}{\bar{\delta}} \beta (t - t^*)$.
One of the roots has been discarded, since it does not satisfy  $p^*(t)\geq \alpha$ for all $t$. The value of $p^*(t)$ for $t\in(t^*,N/s]$ is obtained analogously. 

We notice that  $\min_t p^*(t)=\alpha$ and it is reached for $t=t^*$. This value cannot  further decrease, since $p^*(t)\geq \alpha$.  Thus, from definition in~\eqref{eq:p_11}, different feasible choices of  $\tilde g_1>g_1$ and $\tilde g_2>g_2$ lead to $\tilde p(t) > p^*(t)$. Therefore, the necessary monetary budget $M_\$$ to implement $\tilde p$ instead of  $ p^*$ increases.

Recalling  that $t^*=\gamma N / [s(\beta+\gamma)]$ (see \cite[eq.~11]{sun2020managing}), that $\psi(t)$ is   as  in \eqref{eq:incentive_p1}, and that $r(t)=s$,   we derive the perceived budget from~\eqref{eq:M_per} as follows  
\begin{align}
\label{eq:M_p_star}
    M_{\tup {per}}^* = \frac{\beta\gamma N^2}{2s(\beta+\gamma)} .
\end{align}
Then the final value of the necessary monetary budget $M_\$^*$ follows from~\eqref{eq:Delta_M1}.\hfill\QED

\begin{figure*}[h]\hrule
\begin{equation}
    \label{eq:M_eu_and_M_p_final}
    \begin{split}
    M_\$ & = M_{\tup{per}} + \frac{s}{2}\left[\frac{\bar\delta}{\beta}\int^{\alpha+\frac{t^\ell\beta}{\bar\delta}}_{\alpha}2 \sqrt{y^2-\alpha^2} \,dy + \frac{\gamma}{\bar \delta}\int^{\alpha+\frac{\gamma}{\bar\delta}\left(\frac{N}{s}-t^r\right)}_{\alpha} 2\sqrt{z^2-\alpha^2} \,dz  \right]\\
    &= M_{\tup{per}} + \frac{s\bar\delta}{2\beta}\left.\left(
y\sqrt{y^2-\alpha^2}-\alpha^2\log\left( \sqrt{y^2-\alpha^2} +y \right)    
    \right)\right|^{\alpha+\frac{t^\ell\beta}{\bar\delta}}_{\alpha}\\
    &\quad + \frac{s\bar\delta}{2\gamma}\left.\left(
z\sqrt{z^2-\alpha^2}-\alpha^2\log\left( \sqrt{z^2-\alpha^2} +z \right)    
    \right)\right|^{\alpha+\frac{\gamma}{\bar\delta}\left(\frac{N}{s}-t^r\right)}_{\alpha} 
    \end{split}
\end{equation}
\hrule
\end{figure*}
\subsection{Derivation of $M_\$ = f(M_{\tup{per}})$}\label{app:deriv_f}

Using \eqref{eq:Delta_M1}, we can express $M_\$$ in terms of $M_{\tup {per}}$ as follows:
\begin{equation}
    \label{eq:M_eu_and_M_p_1}
    \begin{split}
    M_\$ &= M_{\tup {per}} + \Delta_M\\
    &=M_{\tup {per}} + \int_0^{\frac{N}{s}}r(\tau)\frac{p(\tau)^2-\alpha^2}{2p(\tau)} \,d\tau .
    \end{split}
\end{equation}

Here, we notice that the interval $[0, \frac{N}{s}]$ can be divided into three sub-intervals using  \eqref{eq:p_limited_budget}. For $t \in (t^l, t^r)$, the discount is $p(t) = 0$ and consequently  the integral is  $0$. Before $t^\ell$ and after $t^r$, the rate $r(t)=s$, as derived in Proposition \ref{prop:opt_disc_lb} and given in \eqref{eq:departure_rate}. Using this, we obtain the following:
\smallskip    
\begin{equation} \label{eq:appendix_definition_f}
M_\$=M_{\tup {per}} + s\int_0^{t^l}\frac{p(\tau)^2-\alpha^2}{2p(\tau)}\,d\tau+s\int_{t^r}^{\frac{N}{s}}\frac{p(\tau)^2-\alpha^2}{2p(\tau)}\,d\tau.
\end{equation}

This equation can be further rewritten by substituting the optimal discount as defined in~\eqref{eq:p_limited_budget}. Furthermore, by performing the  change of variables $y\coloneqq \alpha-\frac{\beta}{\bar\delta} (\tau-t^\ell) $ and $z\coloneqq \frac{\gamma}{\bar\delta} (\tau-t^r)+\alpha $, we obtain \eqref{eq:M_eu_and_M_p_final}. 
Since no real simplification can be performed, we omit the final formulation attained after the  final substitution of by $t^l$ and $t^r$, as in  \eqref{eq:t_l_and_t_r}, where $f(M_\$)$ is replaced by  $M_{\tup{per}}$. Henceforth, we use~\eqref{eq:M_eu_and_M_p_final} to define $M_\$ = f(M_{\tup {per}})$.

\subsection{Proof of Proposition~\ref{prop:opt_disc_lb}}\label{subsec:appendix_proof_2}
The proof  follows similar steps as that for Proposition~\ref{prop:opt_policy}. The lowest travel cost $ \hat C_{\min}\coloneqq \min_t  \hat C(t)$, is experienced by those that receive the lowest incentive, i.e., $0$ if the incentive is the one minimizing the used budget. Following similar steps to~\cite[Sec.~3.2]{sun2020managing}, we notice that, in order for the system to be at optimum, all the commuters have to experience the same total cost $ C(t)= C_e$ for all $t\in[0,N/s]$, this entails that $C(t)=\hat C(t) + C_{\tup{ch}}(t) = C_e =  C_{\min}$, see \eqref{eq:cost_charge}. Therefore, 
\smallskip
\begin{equation}
\label{eq:phi_e_M_p}
\begin{split}
        \int_0^{\frac{N}{s}} r(\tau) \hat C(\tau)  \,d\tau  & +\int_0^{\frac{N}{s}} r(\tau)\hat C_{\tup{ch}}(\tau) \,d\tau  = NC_e \\
        \int_0^{\frac{N}{s}} r(\tau)\hat C(\tau)\,d\tau &= M_{\tup {per}}+N(C_e-\bar \delta \bar p).
\end{split}
\end{equation}
The relation above implies that minimizing  the total travel cost, given $M_{\tup {per}}<M^*_{\tup {per}}$, corresponds to minimizing $C_e$.
By using \eqref{eq:Delta_M1}, we recast the above relation as  follows:
\begin{equation}
    \label{eq:C_e_M_eu_lb}
            \int_0^{\frac{N}{s}} r(\tau)\hat C(\tau)\,d\tau +\Delta_M = M_{\$}+N(C_e-\bar \delta \bar p).
\end{equation}
So, minimizing $C_e$, given a $M_\$<M^*_\$$, coincides with minimizing simultaneously the total travel cost and the inefficiency gap. From the results in \cite[App.~D]{sun2020managing}, it follows that the departure rate that minimises $C_e$ takes the form
\begin{equation} \label{eq:departure_rate}
    r(t)=\begin{cases}
    s & \text{for } t\in[0,t^l]\cup [t^r,N/s]\\
    \alpha s/(\alpha-\beta) & \text{for } t\in(t^l,t'']\\
    \alpha s/(\alpha+\gamma) & \text{for } t\in(t'',t^r),
    \end{cases}
\end{equation}
where $t^l$ and $t^r$ are the starting and ending point of the congestion respectively and $t^l<t''<t^r$, where $t''$ is the time entering at which the bottleneck must be entered to arrive at $t^*$. 
Consequently, the optimal incentive becomes
\begin{equation}\label{eq:psi_lim_budget}
\psi(t)=\begin{cases}    \beta (t^*-t)- C_e & \text{for } t [0,t^l]\\
    0 & \text{for } t\in(t^l,t^r)\\
    \gamma (t-t^*)-C_e & \text{for } t\in[t^r,N/s].
\end{cases}
\end{equation}
From the relation above and the continuity of $\psi(t)$, see \cite[App.~B-C]{sun2020managing}, one can notice that $\varphi_e = \beta (t^*-t^l)= \gamma (t^r-t^*)$. Combining this with the definition of $M_{\tup {per}}$ in \eqref{eq:M_p}, where $\hat \psi(t)$ is as in~\eqref{eq:psi_lim_budget} allows us to compute the analytical expressions for the instants at which the congestion starts $t^\ell$ and ends $t^r$,  and for $t''$. All these expression are derived as functions of $M_{\tup{per}}$. To achieve the final formulation provided in the  proposition's statement, we use the inverse of $f$, computed in~\ref{app:deriv_f}, to obtain these quantities as functions of the monetary budget, i.e., $M_{\tup{per}} = f^{-1}(M_\$)$.

To complete the proof we have to extract from the optimal policy $\psi(t)$ the associated price discount $p(t)$ that is actually used as incentive. As done in Appendix~\ref{app:appendix_proof_1}, we derive it by solving the equation $$C_{\tup{ch}}(t)-\bar\delta\bar p = -\psi(t),$$
where $\psi(t)$ is as in~\eqref{eq:psi_lim_budget}. The resulting final policy $p(t)$, reported in the proposition statement, is guaranteed to minimise the total travel cost of the users and the inefficiency gap, this entails its optimality. \hfill\QED

\addtolength{\textheight}{-12cm}   







\bibliography{references.bib}
\bibliographystyle{ieeetr}

\end{document}